\documentclass[12pt]{article}
\usepackage{amssymb}

\usepackage{amsthm}
\usepackage{amsmath}
\usepackage{graphicx}

 \newtheorem{thm}{Theorem}[section]

 \theoremstyle{definition}

 \numberwithin{equation}{section}
 \def\bR{\mathbb{R}}
\def\bS{\mathbb{S}}
\def\bN{\mathbb{N}}

\def\bu{\mathbf{u}}
\def\bd{\mathbf{d}}

\newtheorem{theorem}{Theorem}[section]

\theoremstyle{definition}

\theoremstyle{remark}
\newtheorem{remark}[theorem]{Remark}

\begin{document}
\title{On a Family of Exact Solutions to the Incompressible Liquid Crystals in Two Dimensions}
\author{Hongjie Dong\footnote{Division of Applied Mathematics, Brown University,
182 George Street, Providence, RI 02912, USA. {\it Email: Hongjie\_Dong@brown.edu}}
\and Zhen
 Lei\footnote{School of Mathematical Sciences; LMNS and Shanghai
 Key Laboratory for Contemporary Applied Mathematics, Fudan University, Shanghai 200433, P. R.China. {\it Email:
 leizhn@gmail.com}}}
\date{\today}
\maketitle

\begin{abstract}
In this paper we construct a family of exact strong
solutions to the
two-dimensional incompressible liquid crystal equations with
finite energy. The initial velocity is chosen to be rotationally
symmetric and the image of the initial orientation of the liquid crystal is a
non-trivial curve on the unit sphere. It turns out that this
family of initial data evolves globally in time by liquid crystal
flow and may shrink to a single point as time goes to infinity.
\end{abstract}

\maketitle





\section{Introduction}

We consider the following hydrodynamic system modelling the flow
of liquid crystal materials in two dimensions (see, for instance,
\cite{Erickson62, Leslie, LiuLin95}):
\begin{equation}\label{LC}
\begin{cases}
\bu_t + \bu\cdot\nabla \bu + \nabla p = \Delta \bu - \nabla\cdot(\nabla
\bd\otimes \nabla \bd),\\[-4mm]\\
\bd_t + \bu\cdot\nabla \bd = \Delta \bd + |\nabla \bd|^2\bd,\\[-4mm]\\
\nabla\cdot \bu = 0,\quad |\bd| = 1,
\end{cases}
\end{equation}
where $\bu$ is the velocity field, $p$ is the scalar pressure and
$\bd$ is the unit-vector on the sphere $\mathbb{S}^2 \subset
\mathbb{R}^3$ representing the macroscopic molecular orientation
of the liquid crystal materials. Here the $i$th component of
$\nabla\cdot(\nabla \bd\otimes \nabla \bd)$ is given by
$\nabla_j(\nabla_i \bd\cdot \nabla_j\bd)$. For simplicity, we have
set all the positive constants in the system to be one. We are
interested in the Cauchy problem of \eqref{LC} with the initial
data
\begin{equation*}
\bu(0, x) = \bu_0(x),\quad \bd(0, x) = \bd_0(x).
\end{equation*}

The above system \eqref{LC} is a simplified version of the
Ericksen--Leslie model for the hydrodynamics of nematic liquid
crystals \cite{Erickson62, Leslie}. The mathematical analysis of
the liquid crystal flow was initiated by Lin and Liu in
\cite{LiuLin95, LiuLin2}. The existence of weak solutions in two
dimensions is obtained in \cite{LinLinWang} where the authors also
showed that there are at most finitely many time singularities for
their weak solutions (see also \cite{Hong}). The uniqueness of
weak solutions in two dimensions was studied in \cite{LinWang,
Zhang}. In \cite{Wang}, the global existence of strong solutions
is proved for initial data with sufficiently small norms of
$\bu_0$ and $\nabla \bd_0$ in ${\rm BMO}^{-1}$. See also
\cite{LiWang} for a small data global existence result in 3D. The
global existence of large solutions is obtained in two dimensions
in \cite{DingLin, LinLinWang} under the assumption that the image
of $\bd_0$ is contained in a half sphere (see also a recent new
proof in \cite{LeiLiZhang}).

In this paper, we are concerned with a family of exact strong solutions with large initial data, which are global in time. We succeed in
constructing them by choosing the rotationally symmetric initial
velocity fields and symmetric initial orientations. Thus the image
of the orientation must be a curve on $\mathbb{S}^2$ which could
be non-trivial. We also show that the image curve of $\bd$ on $\mathbb{S}^2$ may
finally shrink to a single point as time goes to infinity.

The main result of the paper is the following Theorem:
\begin{thm}\label{thm}
Let $\delta_1\in (0,\pi/2)$  and
$u_0, \psi_0$ be two given functions of $r\in \overline{\bR^+}$ which satisfy
\begin{equation}\label{data-p}
\int_0^\infty u_0^2(r)r\,dr<\infty,\quad \int_0^\infty (\psi_0'(r))^2r\,dr<\infty
\quad
\delta_1 < \psi_0(r) <\pi -\delta_1.
\end{equation}
Then there exists a smooth function $\Phi\,:\,(\delta_1,\pi-\delta_1)\to \bR$ such that, for the rotationally symmetric initial velocity field
\begin{equation}\label{constraint-u}
\bu_0(x) = u_0(|x|)e_\theta,\quad e_\theta = (- \frac{x_2}{|x|},
\frac{x_1}{|x|})^T
\end{equation}
and the symmetric initial orientation unit vector
\begin{equation}\label{constraint-d}
\bd_0(x) = \begin{pmatrix}
             \sin\psi_0(|x|)\cos\phi_0(|x|) \\
             \sin\psi_0(|x|)\sin\phi_0(|x|) \\
             \cos\psi_0(|x|)\\
           \end{pmatrix}, \quad \phi_0 = \Phi(\psi_0),
\end{equation}
the solution to the liquid crystal equations \eqref{LC} with
initial data \eqref{constraint-u}-\eqref{constraint-d} is global
and can be exactly solved.
\end{thm}

\begin{remark}
The solution we obtained in Theorem \ref{thm} becomes smooth instantaneously for $t>0$. By the weak-strong uniqueness (see, for instance, \cite{LinWang} and \cite{Zhang}), such solution is also unique.
\end{remark}

We point out here that the $L^2$ norms of $\bu_0$ and $\nabla \bd_0$ may
be large in our theorem. Hence the global well-posedness of such
solutions is not a consequence of \cite{Wang}. Since the image of
$\bd_0$ is not contained in any half sphere on $\mathbb{S}^2$ (cf. \eqref{eq12.26}), the
existence of our solutions cannot be deduced from
\cite{DingLin, LinLinWang, LeiLiZhang} either. The main ingredient of this paper is
that we can give an explicit expression of the solution to the
liquid crystal equations \eqref{LC}  with
initial data \eqref{constraint-u}-\eqref{constraint-d}. For the proof, we first
parameterize the curve, and solve a nonlinear ordinary
differential equation \eqref{eq10.47} to reduce the problem to a nonlinear heat equation. We then construct a family of exact solutions by solving the nonlinear heat equation using a Hopf--Cole type transformation. We also show that the image curve of $\bd$ has non-increasing arc length along the liquid crystal flow. However, the arc length does not always decreasing in general. We will give an example of global solutions which does not shrink in time at the end of this paper.

Before ending this introduction, let us mention a related result
on harmonic map heat flow in \cite{ChangDingYe} where finite-time
singularities are shown for a class of initial data in the form
\begin{equation}\nonumber
\bd_0(x) = \begin{pmatrix}
           x_1r^{-1}\sin\psi_0(r) \\
           x_2r^{-1}\sin\psi_0(r) \\
           \cos\psi_0(r) \\
         \end{pmatrix}
\end{equation}
with
\begin{equation}\nonumber
\psi_0(0) = 0,\ \psi_0(R) > \pi\ {\rm for\ some}\  R > 0.
\end{equation}
This also serves as an example of finite-time singularities for the
liquid crystal system \eqref{LC} in the special case when $\bu\equiv 0$.  Whether there are finite-time singularities
for the incompressible liquid crystal flow in two dimensions with
finite energy and non-trivial velocity remains a very interesting open question.

The paper is organized as follows: In Section 2 we will
parameterize the orientation unit vector $\bd$ and rewrite the
liquid crystal equations \eqref{LC} in terms of $\phi$ and $\psi$.
Then in Section 3 we solve \eqref{LC} by giving an exactly
expression of $\bu$ and $\bd$, and finish the proof of Theorem
\ref{thm}.

\section{Parameterization of the Liquid Crystal Equations}

In this section we parameterize the liquid crystal equations. Let
$e_r=x/|x|$. It is clear that the gradient operator can be
expressed as
\begin{equation}\nonumber
\nabla = e_r\partial_r + r^{-1}e_\theta\partial_\theta,\quad
\partial_\theta = x_1\partial_2 - x_2\partial_1.
\end{equation}
We look for solution of $(\bu, \bd, p)$, where $\bu$ is rotational
and $p, \bd$ are symmetric, i.e.,
$$
\bd=\bd(t,r),\quad \bu=u(t,r)e_\theta,\quad p=p(t,r).
$$
We can rewrite  \eqref{LC} into the following system
\begin{equation}
                            \label{eq9.56}
\begin{cases}
u_t=u_{rr}+\frac 1 r u_r-\frac 1 {r^2} u \\[-4mm]\\
    p_r=\frac 1 r u^2-\frac 1 2 (|\bd_r|^2)_r-\Delta \bd\cdot \bd_r\\[-4mm]\\
    \bd_t=\bd_{rr}+\frac 1 r \bd_r+|\bd_r|^2 \bd
\end{cases}
\end{equation}
Since $\bd\in \mathbb{S}^2$, we may assume
\begin{equation}
                                \label{eq10.01}
\bd=\begin{pmatrix}
    \sin\psi\cos\phi \\
   \sin\psi\sin\phi \\
    \cos\psi \\
  \end{pmatrix}.
\end{equation}
A simple computation gives
\begin{equation}
                                    \label{eq10.02}
\bd_\psi=\begin{pmatrix}
           \cos\psi\cos\phi \\
           \cos\psi\sin\phi \\
           -\sin\psi \\
         \end{pmatrix},
\quad \bd_\phi=\begin{pmatrix}
        -\sin\psi\sin\phi \\
        \sin\psi\cos\phi \\
        0 \\
      \end{pmatrix},
\end{equation}
and
\begin{equation}
                                        \label{eq10.04}
\bd_{\psi\psi}=-\bd, \quad \bd_{\psi\phi}=
\begin{pmatrix}
         -\cos\psi\sin\phi \\
         \cos\psi\cos\phi \\
        0\\
\end{pmatrix},\quad
\bd_{\phi\phi}=\begin{pmatrix}
              -\sin\psi\cos\phi\\
              -\sin\psi\sin\phi \\
              0 \\
            \end{pmatrix}.
\end{equation}
Clearly, from \eqref{eq10.01}, \eqref{eq10.02}, and
\eqref{eq10.04}, we have
\begin{equation}
                            \label{eq10.18}
\begin{split}
  \bd_\psi\cdot \bd=0,\quad \bd_\psi\cdot \bd_\psi=1,\quad \bd_\psi\cdot \bd_\phi=0,\\
  \bd_\psi\cdot \bd_{\psi\psi}=0,\quad \bd_{\psi}\cdot
\bd_{\psi\phi}=0,\quad \bd_{\psi}\cdot
\bd_{\phi\phi}=-\cos\psi\sin\psi,
\end{split}
\end{equation}
and
\begin{equation}
                            \label{eq10.20}
\begin{split}
 \bd_\phi\cdot \bd=0,\quad \bd_\phi\cdot \bd_\psi=0,\quad \bd_\phi\cdot \bd_\phi=\sin^2\psi, \\
\bd_\phi\cdot \bd_{\psi\psi}=0,\quad \bd_{\phi}\cdot
\bd_{\psi\phi}=\sin\psi\cos\psi, \quad \bd_{\phi}\cdot
\bd_{\phi\phi}=0.
\end{split}
\end{equation}
By the chain rule,
\begin{align*}
\bd_t&=\bd_\psi\psi_t+\bd_\phi\phi_t,\\
\bd_r&=\bd_\psi\psi_r+\bd_\phi\phi_r,\\
\bd_{rr}&=\bd_{\psi\psi}(\psi_r)^2+2\bd_{\psi\phi}\psi_r\phi_r+\bd_{\phi\phi}(\phi_r)^2
+\bd_\psi\psi_{rr}+\bd_\phi\phi_{rr}.
\end{align*}
Combining the above equalities with the third equation of
\eqref{eq9.56}, we obtain
\begin{align}
\bd_\psi\psi_t+\bd_\phi\phi_t=\bd_{\psi\psi}(\psi_r)^2+2\bd_{\psi\phi}\psi_r\phi_r+\bd_{\phi\phi}(\phi_r)^2
+\bd_\psi\psi_{rr}\nonumber\\
                             \label{eq10.29}
 +\bd_\phi\phi_{rr}+\frac 1 r
(\bd_\psi\psi_r+\bd_\phi\phi_r)+\bd(\psi_r^2+\sin^2\psi\phi_r^2).
\end{align}
We dot \eqref{eq10.29} with $\bd_\psi$ and $\bd_\phi$ respectively
and use \eqref{eq10.18} and \eqref{eq10.20} to deduce the
following two equations
\begin{align}
                                            \label{eq10.35}
\psi_t&=\psi_{rr}+\frac 1 r \psi_r-\cos\psi\sin\psi\phi_r^2,\\
                                        \label{eq10.37}
\phi_t\sin^2 \psi&=(\phi_{rr}+\frac 1 r \phi_r)\sin^2
\psi+2\cos\psi \sin\psi\psi_r\phi_r.
\end{align}
Note that in \eqref{eq10.29} the terms in the $\bd$-direction cancel
each other.

We assume that the image of $\bd$ is a curve on $\bS^2$. In other words, $\phi$ and $\psi$ are not independent.
Take $\phi=\Phi(\psi)$. Then it follows from
\eqref{eq10.35} and \eqref{eq10.37} that
\begin{equation}
                                    \label{eq11.10}
\psi_t=\psi_{rr}+\frac 1 r
\psi_r-(\Phi')^2\cos\psi\sin\psi\psi_r^2,
\end{equation}
$$
\Phi'\psi_t\sin^2 \psi=\Phi'(\psi_{rr}+\frac 1 r
\psi_r)\sin^2 \psi+\Phi''(\psi_r)^2\sin^2 \psi+
2\Phi'\cos\psi \sin\psi(\psi_r)^2.
$$
By substituting the expression of $\psi_t$ in \eqref{eq11.10}
into the last equation, we then deduce
\begin{equation}
                                        \label{eq10.44}
(\Phi')^3\cos\psi\sin^3\psi(\psi_r)^2+\Phi''\sin^2\psi(\psi_r)^2
+2\Phi'\cos\psi \sin\psi(\psi_r)^2=0.
\end{equation}
A sufficient condition for \eqref{eq10.44} is the following ODE
\begin{equation}
                                    \label{eq10.47}
\Phi''\sin^2\psi+2 \Phi'\cos\psi
\sin\psi+(\Phi')^3\cos\psi\sin^3\psi=0.
\end{equation}

\section{Proof of Theorem \ref{thm}}
In this section, we complete the proof of Theorem \ref{thm} by explicitly solving the nonlinear ODE \eqref{eq10.47} and a nonlinear heat equation.
Note that \eqref{eq10.47} is equivalent to
$$
(\sin^2\psi\Phi')'+\cos\psi\sin^3\psi (\Phi')^3=0.
$$
Thus, for $\psi\neq 0,\pi$, we have
$$
\frac {(\sin^2\psi\Phi')'}{(\sin^2\psi\Phi')^3}=-\frac
{\cos\psi}{\sin^3\psi},
$$
and then
$$
\frac 1 {(\sin^2\psi \Phi')^2}=-\frac 1 {\sin^2\psi}+\beta
$$
for some constant $\beta>1$. Consequently,
\begin{equation}
                                    \label{eq11.50}
(\Phi')^2=\frac 1 {(\beta\sin^2\psi-1)\sin^2\psi}.
\end{equation}
Here we require
\begin{equation}
                                    \label{eq12.05}
\beta\sin^2\psi> 1.
\end{equation}
Recall the
condition on the initial data
\begin{equation}
                            \label{eq11.57}
\psi_0\in (\delta_1,\pi-\delta_1).
\end{equation}
We fix a $\beta\ge 1/\sin^2\delta_1$.
Plugging the expression \eqref{eq11.50} into \eqref{eq11.10}, we
get
\begin{equation}
                            \label{eq11.16}
\psi_t=\psi_{rr}+\frac 1 r \psi_r-\frac {\cos\psi}
{(\beta\sin^2\psi-1)\sin\psi}\psi_r^2,
\end{equation}
which is a nonlinear heat equation in 2D. This motivates us to use the
idea of the Hopf--Cole transformation to reduce \eqref{eq11.16} to a linear
equation. Let $F=F(\psi)$ be a function to be chosen. Clearly,
$$
F_t=F'\psi_t,\quad F_{rr}+\frac 1 r F_r=F'(\psi_{rr}+\frac 1 r
\psi_r)+F''\psi_r^2.
$$
Then $F$ satisfies the linear heat equation
\begin{equation}
                                    \label{eq11.28}
F_t=F_{rr}+\frac 1 r F_r
\end{equation}
provided that
$$
F''+F'\frac {\cos\psi} {(\beta\sin^2\psi-1)\sin\psi}=0.
$$
We solve the above ODE to get
$$
F'(\psi)=C_1\frac{\sin\psi}{(\beta\sin^2\psi-1)^{1/2}},
$$
and
\begin{equation*}
F(\psi)=\int
C_1\frac{\sin\psi}{(\beta\sin^2\psi-1)^{1/2}}\,d\psi =\frac
1 {\beta^{1/2}}\left(C_1 \arccos\Big(\sqrt{\frac \beta
{\beta-1}}\cos\psi\Big)+C_2\right).
\end{equation*}
Here $C_1$ and $C_2$ are arbitrary constants. For our purpose, it
is convenient to take $C_1=\beta^{1/2}$ and $C_2=0$. Then
$$
F(\psi)=\arccos\Big(\sqrt{\frac \beta {\beta-1}}\cos\psi\Big).
$$
Note that $F$ is well defined for any $\psi\in (\delta_1,\pi-\delta_1)$,
and it is a strictly increasing continuous function in the same interval.
By \eqref{eq11.57}, we have
$$
F(\psi_0)\in (\delta_2,\pi-\delta_2),\quad \delta_2=\arccos\Big(\sqrt{\frac \beta {\beta-1}}\cos\delta_1\Big)\in (0,\pi/2).
$$
From \eqref{eq11.28}, we have
$$
F(\psi(t,\cdot))=\Gamma(t,\cdot)*F(\psi_0),
$$
where $\Gamma(t,\cdot)$ is the 2D heat kernel.
Hence,
$$
F(\psi(t,\cdot))\in (\delta_2,\pi-\delta_2),\quad \psi(t,\cdot)\in (\delta_1,\pi-\delta_1),
$$
which implies that the condition \eqref{eq12.05} is satisfied for any $t\ge 0$ if it is satisfied at $t=0$. Therefore, the solution $\psi$ exists globally. Furthermore, by the maximum principle, the length of the image curve of $d$ is non-increasing.

It follows from \eqref{eq11.50} that
$$
\Phi'(\psi)=\frac 1 {(\beta \sin^2 \psi-1)^{1/2} \sin\psi}.
$$
Here we took the positive branch of the square root. Therefore, by choosing a suitable constant,
$$
\Phi(\psi)=\int_{\pi/2}^\psi\frac 1 {(\beta \sin^2 s-1)^{1/2} \sin s}\,ds.
$$
Observe that
\begin{equation}
                                    \label{eq12.26}
\lim_{\psi\to \delta_1}\Phi(\psi)\to -\infty,\quad
\lim_{\psi\to \pi-\delta_1}\Phi(\psi)\to +\infty\quad \text{as}\,\,\delta_1\to 0.
\end{equation}

We now treat the first two equations in \eqref{eq9.56}. In the equation of $u$, we divide both sides by $r$ and get
\begin{equation}
                            \label{eq2.56}
\left(\frac u r\right)_t=\left(\frac u r\right)_{rr}+\frac 3 r \left(\frac u r\right)_r.
\end{equation}
We note that $\partial_{rr}+\frac 3 r\partial_r$ is the Laplace operator in $\bR^4$ for radially symmetric functions. Moreover, as a radial function in $\bR^4$, $u_0/r\in L^2(\bR^4)$ since, by \eqref{data-p},
$$
\int_{\bR^4}u_0^2(|y|)|y|^{-2}\,dy=C\int_{0}^\infty u_0^2(r)r\,dr<\infty.
$$
Therefore, due to \eqref{eq2.56} we get
$$
 u (t,r)=r\int_{\bR^4}\frac 1 {(4\pi t)^2}e^{-\frac {|x-y|^2}{4t}}u_0 (t,|y|)|y|^{-1}\,dy,\quad x=r(1,0,0,0),
$$
and, as a radial function in $\bR^4$, $u(t,\cdot)/r\in L^2(\bR^4)$ for any $t\ge 0$. This further implies that $\bu(t,\cdot)\in L^2(\bR^2)$. Since $\bu=(-x_2,x_1)u/r$, we see that $\bu$ is a smooth function. Finally, from the second equation of \eqref{eq9.56}, we get
$$
p_r=\frac 1 r u^2-\bd_r\cdot \bd_{rr}-(\bd_{rr}+\frac 1 r \bd_r)\cdot \bd_r=\frac 1 r u^2-2\bd_r\cdot \bd_{rr}-\frac 1 r \bd_r\cdot \bd_r.
$$
Therefore,
$$
p(t,r)=\int_0^r \left(\frac 1 s u^2(t,s)-\frac 1 s |\bd_r(t,s)|^2\right)\,ds-|\bd_r(t,r)|^2.
$$
Note that by the radial symmetry condition, the integral above is convergence. This completes the proof of Theorem \ref{thm}.

\begin{remark}
In the special case that the limit $\psi_0$ exists as $r\to \infty$, from the proof above it is easily seen that $\psi(t,\cdot)$ converges to this limit uniformly in $\bR^2$ as $t$ goes to infinity.
Indeed, $F(\psi_0)$ is a bounded function and has a limit $\bar F$ as $r\to \infty$. By the simple property of solutions to the heat equation, we have
$$
\lim_{t\to \infty}F(\psi(t,\cdot))\to \bar F
$$
uniformly in $\bR^2$, which implies that
$$
\lim_{t\to \infty}\psi(t,\cdot)\to \bar\psi:=\arccos\Big(\sqrt{\frac {\beta-1} \beta }\cos\bar F\Big)
$$
uniformly in $\bR^2$. In this case, the image curve shrinks to a point on $\bS^2$ uniformly along the liquid crystal flow.

On the other hand, one can find certain initial data such that the image curve on $\bS^2$ does not shrink at all.
This can be seen from the one-to-one correspondence of $F$ and $\psi$ and the following example. Let $v_0$ be a radial function in $\bR^2$ defined as
$$
v_0(r)=0\,\, \text{on}\,\,[0,e],\quad
v'_0(r)=\begin{cases} \frac 1 r\frac 1 {(6k+2)^3-(6k+1)^3}\quad \text{on}\,\,(e^{(6k+1)^3},e^{(6k+2)^3}]\\[-4mm]\\
0\quad \text{on}\,\,(e^{(6k+2)^3},e^{(6k+4)^3}]\cup (e^{(6k+5)^3},e^{(6k+7)^3}]\\[-4mm]\\
-\frac 1 r\frac 1 {(6k+5)^3-(6k+4)^3}\quad \text{on}\,\,(e^{(6k+4)^3},e^{(6k+5)^3}]
\end{cases},\,\,k\in \bN.
$$
It is easy to check that $\nabla v_0\in L^2(\bR^2)$, $v_0\in [0,1]$, and
$$
v_0(r)=
\begin{cases} 1\quad \text{on}\,\,[e^{(6k+2)^3},e^{(6k+4)^3}]\\[-4mm]\\
0\quad \text{on}\,\,[e^{(6k+5)^3},e^{(6k+7)^3}]
\end{cases},\,\,k\in \bN.
$$
Now take $t_k=e^{2(6k+3)^3}$, $\tilde t_k=e^{2(6k+6)^3}$,  and let $v$ be the solution to the heat equation with initial data $v_0$. Since
$$
\lim_{k\to \infty}\frac {e^{2(6k+2)^3}}{t_k}=0, \quad
\lim_{k\to \infty}\frac {e^{2(6k+4)^3}}{t_k}=\infty,
$$
we get
$$
\lim_{k\to \infty} \int_{|x|\in (e^{(6k+2)^3},e^{(6k+4)^3})}\Gamma(t_k,x)\,dx=1.
$$
Therefore,
$$
\lim_{k\to \infty}v(t_k,0)=1.
$$
Similarly, we have
$$
\lim_{k\to \infty}v(\tilde t_k,0)=0.
$$
Thus, the range of $v$ does not shrink in time.
\end{remark}

\section*{Acknowledgement}
This work was done when Zhen Lei was visiting the Division of Applied
Mathematics of Brown University during 2012. He would like to
thank the hospitality of the institute. Hongjie Dong was partially supported by the NSF under agreement DMS-0800129 and DMS-1056737. Zhen Lei was
supported by NSFC (grant No.11171072), the Foundation for
Innovative Research Groups of NSFC (grant No.11121101), FANEDD,
Innovation Program of Shanghai Municipal Education Commission (grant
No.12ZZ012),  NTTBBRT of NSF (No. J1103105), and SGST 09DZ2272900.


\frenchspacing
\bibliographystyle{plain}

\begin{thebibliography}{99}

\bibitem{ChangDingYe} K. Chang, W. Ding and R. Ye, \textit{Finite-time blow-up of the
heat flow of harmonic maps from surfaces}. J. Differ. Geom.
\textbf{36} (1992), 507--515.

\bibitem{DingLin} W. Ding and F.-H. Lin, \textit{A generalization of Eells-Sampson's
theorem}. J. Partial Differential Equations 5 (1992), no. 4,
13--22.



\bibitem{Erickson62} J. Ericksen, \textit{Hydrostatic theory of liquid crystal}.
Arch. Ration. Mech. Anal. \textbf{9 }(1962), 371--378.

\bibitem{Hong} Min-Chun Hong,  \textit{Global existence of solutions of the simplified
Ericksen-Leslie system in dimension two}. Calc. Var. Partial
Differential Equations \textbf{40} (2011), no. 1-2, 15--36.

\bibitem{LeiLiZhang} Z. Lei, D. Li and X. Zhang, \textit{A new proof of global wellposedness of liquid
crystals and heat harmonic maps in two dimensions}. preprint.

\bibitem{Leslie} F. Leslie, \textit{Some constitutive equations for liquid crystals}.
Arch. Rational Mech. Anal. \textbf{28} (1968), no. 4, 265--283.

\bibitem{LiuLin95} F.-H. Lin and C. Liu, \textit{Nonparabolic dissipative systems modeling the flow of liquid crystals}. Comm. Pure Appl. Math. 48 (1995), no. 5, 501--537.

\bibitem{LiuLin2} F.-H. Lin and C. Liu, \textit{Partial regularity of the dynamic system modeling the flow of liquid crystals}. Discrete Contin. Dynam. Systems 2 (1996), no. 1, 1--22.

\bibitem{LiWang} X. Li and D. Wang,  \textit{Global solution to the incompressible flow
of liquid crystals}. J. Differential Equations \textbf{252}
(2012), no. 1, 745--767.

\bibitem{LinLinWang} F.-H. Lin, J. Lin and C. Wang, \textit{Liquid crystal flows in two dimensions}. Arch. Ration. Mech. Anal. \textbf{197} (2010),
no. 1, 297--336.

\bibitem{LinWang} F.-H. Lin and C. Wang,  \textit{On the uniqueness of heat flow of
harmonic maps and hydrodynamic flow of nematic liquid crystals}.
Chin. Ann. Math. Ser. B \textbf{31} (2010), no. 6, 921--938.

\bibitem{Wang} C. Wang, \textit{Well-posedness for the heat flow of harmonic maps and the liquid crystal flow with rough initial data}. Arch. Ration.
Mech. Anal. \textbf{200} (2011), no. 1, 1--19.


\bibitem{Zhang} X. Xu and Z. Zhang,  \textit{Global regularity and uniqueness of weak
solution for the 2-D liquid crystal flows}. J. Differential
Equations \textbf{252} (2012), no. 2, 1169--1181.

\end{thebibliography}

\end{document}